\documentclass[a4paper, 10pt]{article}

\usepackage{natbib}

\usepackage{graphicx}%
\usepackage{multirow}%
\usepackage{amsmath,amssymb,amsfonts}%
\usepackage{amsthm}%
\usepackage{mathrsfs}%
\usepackage[title]{appendix}%
\usepackage{xcolor}%
\usepackage{textcomp}%
\usepackage{manyfoot}%
\usepackage{booktabs}%
\usepackage{algorithm}%
\usepackage{algorithmicx}%
\usepackage{algpseudocode}%
\usepackage{listings}%

\theoremstyle{thmstyleone}%
\newtheorem{theorem}{Theorem}
\newtheorem{proposition}[theorem]{Proposition}%
\newtheorem{lemma}[theorem]{Lemma}%
\newtheorem{corollary}[theorem]{Corollary}%

\theoremstyle{thmstyletwo}%
\newtheorem{remark}{Remark}%

\theoremstyle{thmstylethree}%
\newtheorem{definition}{Definition}%
\numberwithin{equation}{section}

\begin{document}

\title{Scoring Nim}
\author{Hiromi Oginuma\thanks{Graduate School of Humanities and Sciences, Nara Women's University.\\ \qquad xah\_oginuma@cc.nara-wu.ac.jp} \and Masato Shinoda\thanks{Division of Natural Sciences, Nara Women's University.\\
\qquad shinoda@cc.nara-wu.ac.jp}}






\maketitle

\begin{abstract}
Nim is a well-known combinatorial game in which two players alternately remove stones from distinct piles. A player who removes the last stone wins under the normal play rule, while a player loses under the mis\`ere play rule. In this paper, we propose a new variant of Nim with scoring that generalizes both the normal and mis\`ere play versions of Nim as special cases. We study the theoretical aspects of this extended game and analyze its fundamental properties, such as optimal strategies and payoff functions.
\end{abstract}



\section{Introduction}\label{sec1}

Nim is a well-known combinatorial game that has been extensively studied. In the standard version of Nim, there are three piles of stones, and two players take turns. On each turn, a player selects one pile and removes any positive number of stones from that pile. The position of the game is represented as $(x, y, z)$ where $x,y$ and $z$ are non-negative integers representing the number of stones in each pile. 
For example, from $(3, 5, 7)$, the player to move can move to $(3, 2, 7)$, $(0, 5, 7)$, and so on. Players take turns making these moves, and whoever takes the last stone {\it wins} the game, while a player with no available moves automatically loses. This rule is known as the {\it normal} play rule. If the position is $(3, 0, 0)$, the player to move can move to $(0, 0, 0)$, the terminal position, and win the game. This game can be played with any number of piles following the same rules.

In contrast, the rule where the player who takes the last stone {\it loses} is known as the {\it mis\`ere} play rule. Under this rule, if the position is $(3, 0, 0)$, the player to move can move to $(1,0,0)$, thereby forcing the opponent to take the last stone and securing victory.

Nim is a zero-sum, deterministic, perfect-information, impartial and finite game with no possibility of a draw. Therefore, every position can be classified as either an N-position, where the player to move can force a win, or a P-position, where the player to move will inevitably lose if the opponent plays optimally. The classification of N-positions and P-positions in Nim was first established by Bouton (1901/1902)
\cite{bou02}, revealing its mathematical structure. Nim and its variants are discussed in detail for instance by Albert, Nowakowski and Wolfe (2011) \cite{alb11}, Berlekamp, Conway and Guy (2001, 2003, 2004) \cite{ber04} and Siegel (2013) \cite{sie13}.

In this paper, we propose a new variant of Nim named {\it Scoring Nim}, in which the outcome is not determined solely by whether a player takes the last stone, but rather by assigning predetermined bonus points for doing so. In this game, a player earns 1 point per stone taken, and an additional bonus of $N$ points for taking the last stone. The player with the higher total points wins. Here, $N$ is fixed before the game starts, which is allowed to take non-integer values and can be any real number. If $N$ is negative, taking the last stone results in a penalty.
Under this setting, when the bonus $N=\infty$, the game corresponds to the normal play Nim; when $N=0$, it becomes a simple stone-taking game; and when $N=-\infty$, it corresponds to the mis\`ere play Nim. 

Moreover, we assume that the two players not only compete to win but also aim to accumulate as many points as possible. 
We examine the fundamental properties, significance, and the appeal of Scoring Nim. In particular, we analyze how the optimal strategy varies in a complex manner depending on the bonus setting. These results indicate that Scoring Nim possesses intriguing properties worthy of further analysis.

This paper is organized as follows. In Section \ref{sec2} we formally introduce the rules of the new game. In Section \ref{sec3}  we define the payoff function, which plays a crucial role in the analysis of this game. In Section  \ref{sec4}  we examine the optimal strategy for the case of two piles of stones. In Section  \ref{sec5}  we discuss general properties of the payoff function. In Section  \ref{sec6} we fix the number of piles at three and explicitly present the payoff function when the minimum number of stones in a pile is one. We analyze the reasons behind the numerous breakpoints (i.e., points where the increasing and decreasing trends switch) in the plot of the function.

This paper is a compilation of Oginuma and Shinoda (2024) \cite{ogi24a,ogi24b}. 
This resubmission is allowed under the guidelines of the Information Processing Society of Japan 
(For details, see \\{\tt https://www.ipsj.or.jp/english/faq/ronbun-faq\_e.html}).

\section{The rules of Scoring Nim}\label{sec2}

In this section, we provide a detailed explanation of the rules of Scoring Nim.

\medskip

\begin{definition}[Scoring Nim]\label{def1}
Let $N$ be a real number fixed before the game starts. There are $n$ piles of stones, and two players take turns. On each turn, a player selects one pile and takes any positive number of stones from it. This process continues until all stones are taken. Each player earns
\begin{itemize}
\item $1$ point for each stone taken,
\item $N$ points for taking the last stone.
\end{itemize}
The total points accumulated by each player determine the final score.
\end{definition}

\medskip

A game position is represented\footnote{When the number of piles is 1, 2, or 3, the game positions are represented as $(x)$, $(x, y)$, and $(x, y, z)$, respectively.} as ${\boldsymbol p}=(p_{1},p_{2},\ldots,p_{n})$, where each $p_{i}$ is a non-negative integer. The terminal position $(0,0,\ldots,0)$ is denoted by $\boldsymbol{0}$. The total number of stones in a position $\boldsymbol{p}$ is denoted as $|\boldsymbol{p}|=\sum_{i=1}^{n}p_{i}$. As a concrete example, consider a game where the initial position is $(3,2,1)$, and two players (the first and second player) take turns as follows: 
\[(3,2,1)\to (2,2,1)\to (2,2,0)\to(2,0,0)\to(0,0,0).\]
If the game proceeds in this manner, both players take three stones each. The second player takes the last stone and receives the bonus.
Thus, the first player's score is $3$ points, while the second player's score is $3+N$ points.

If the winner is determined by the higher score, then as $N=\infty$, the game corresponds to the normal play Nim, and as $N=-\infty$, it corresponds to the mis\`ere play Nim. In this sense, Scoring Nim serves as an extension of Nim. Moreover, even if the initial position is the same, the player with a winning strategy may differ depending on the value of $N$.

Furthermore, Scoring Nim becomes even more intriguing when the objective shifts from merely securing victory to maximizing points. This setting corresponds to {\it scoring games} in combinatorial game theory (for details, see Stewart (2019) \cite{ste17}). In \cite{ste17}, Section \ref{sec5} primarily discusses the scoring adaptation of octal games as a topic related to Nim\footnote{Additionally, in the introduction of \cite{ste17}, there is a brief mention of a stone-taking game corresponding to Scoring Nim with $N=0$.}. Scoring Nim provides a simple yet fascinating extension of Nim by introducing scores. It serves as a valuable example of scoring games and offers important insights into the study of such games. While our analysis does not explicitly employ the general theoretical framework of scoring play, the game we introduce falls naturally within the class of impartial scoring games considered in the literature, and may provide a useful concrete case for further theoretical investigation.

\begin{figure}[h] 
\centering
\hspace*{.5cm}
\unitlength.15pt
\begin{picture}(1500,500)(0,200)
\put(200,500){\line(5,1){500}}
\put(700,600){\line(5,-1){500}}
\put(700,500){\line(0,1){100}}
\put(600,650){\large{$(5,4,2)$}}
\put(100,430){\large{$(5,4,1)$}}
\put(600,430){\large{$(1,4,2)$}}
\put(1100,430){\large{$(0,4,2)$}}
\put(200,300){\line(0,1){100}}
\put(500,300){\line(2,1){200}}
\put(700,400){\line(1,-1){100}}
\put(1100,300){\line(1,1){100}}
\put(1200,400){\line(2,-1){200}}
\put(100,230){\large{$(2,4,1)$}}
\put(400,230){\large{$(1,3,2)$}}
\put(700,230){\large{$(1,0,2)$}}
\put(1000,230){\large{$(0,2,2)$}}
\put(1300,230){\large{$(0,0,2)$}}
\end{picture}
\caption{A partial game tree from the initial position $(5,4,2)$}
\label{fig1}
\end{figure}
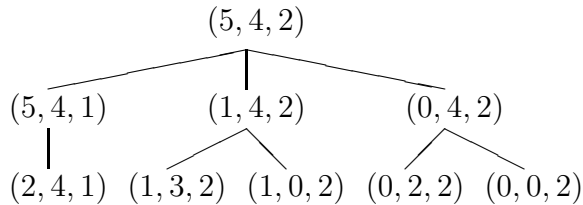

Let us consider Scoring Nim with an initial position of $(5,4,2)$. Players' move selections at each turn can be represented as an extensive-form game tree. By applying the minimax method, one can determine the optimal move for each position and the optimal sequence of moves leading to the terminal position, assuming both players play optimally. Fig. \ref{fig1} presents the game tree for the initial sequence of moves, illustrating some representative promising moves.

The following conclusions are presented here as preliminary observations in order to illustrate the dependence of optimal play on the value of the bonus $N$; a general theoretical framework justifying them will be established in subsequent sections.
In this initial position, the optimal move for the first player is:
\begin{itemize}
\item Moves to $(5,4,1)$ when $N\leq -4$,
\item Moves to $(1,4,2)$ or $(0,4,2)$ when $-4\leq N\leq -1$,
\item Moves to $(0,4,2)$ when $-1\leq N\leq 1$,
\item Moves to $(1,4,2)$ when $1\leq N\leq 4$,
\item Moves to $(5,4,1)$ when $4\leq N$.
\end{itemize}
As shown above, the optimal move selection at each position depends on the value of $N$. When 
$|N|$ is sufficiently large, the score is significantly influenced by whether the last stone is taken. In both the normal and mis\`ere play versions of Nim, the optimal move is to $(5,4,1)$, which we refer to as {\it the Nim-winning strategy}.
On the other hand, when $N=0$, the best move is to maximize the number of stones taken by transitioning to $(0,4,2)$, which we refer to as {\it the greedy strategy}. However, it is important to note that there are also values of $N$ for which other moves become optimal. For example, when $N=3$, the optimal sequence of moves is:
\[(5,4,2)\to (1,4,2)\to (1,3,2)\to(1,0,2)\to(1,0,1)\to(1,0,0)\to(0,0,0).\]
In this case, the first player takes $8$ stones and earns $8$ points, while the second player takes $3$ stones and, with the bonus, obtains $3+N=6$ points. Intuitively, these intermediate strategies can be understood as moves that force the opponent into a position where neither the Nim-winning strategy nor the greedy strategy is particularly effective. A detailed discussion of these intermediate strategies for such values of $N$ will be provided in Section  \ref{sec6}. 

\section{The payoff function and recursion formula}\label{sec3}

Scoring Nim is a constant-sum game in which the total points of both players always add up to $|\boldsymbol{p}|+N$ when starting from the position $\boldsymbol{p}$. Therefore, treating it as a zero-sum game, we assume that each player aims to maximize the final score difference, i.e., 
\[\mbox{Payoff} = (\mbox{First player's score})-(\mbox{Second player's score})\]
at the end of the game. 
The sign of this value indicates whether the player has gained more points than the opponent. 
Let $f_{N}(p_{1},p_{2},\ldots,p_{n})$ denote this score difference for the first player, given in the initial position $\boldsymbol{p}=(p_{1},p_{2},\ldots,p_{n})$. Then, for $n=1$ (i.e., when there is only one pile of stones), we obtain $f_{N}(1)=1+N$ and  
\begin{equation}\label{eq1}
f_{N}(x)=\max\{x+N, x-2-N\}=x-1+|1+N|
\end{equation}
for $x\geq 2$.
From the comparison in (\ref{eq1}), when the initial position is $(x)$ with $x\geq 2$, the first player should take all the stones if $N\geq -1$. On the other hand, if $N\leq -1$, the first player should take $x-1$ stones on the first move, leaving only one stone for the opponent's turn.

Similar to $f_{N}(x)$ with $x\geq 2$, $f_{N}(\boldsymbol{p})$ is recursively defined based on comparisons with subsequent positions. To formalize this, we extend the definition of the payoff function $f_{N}(\boldsymbol{p})$ to positions beyond the initial state. 
We define it as the difference between the total points obtained by the current player and those obtained by the opponent from that position onward.
To ensure consistency, we set $f_{N}(\boldsymbol{0})=-N$, meaning that when no stones remain, the opponent has just received an endgame bonus of $N$ points. 
Let $\boldsymbol{p}\to\boldsymbol{p}^{\prime}$ denote that a transition from position $\boldsymbol{p}$ to $\boldsymbol{p}^{\prime}$ is possible in a single move. 
Then, for a position $\boldsymbol{p}$ that is not a terminal position, the following recurrence relation holds:
\begin{equation} \label{eq2}
f_{N}(\boldsymbol{p})=\max_{\boldsymbol{p}\to\boldsymbol{p}^{\prime}}\{|\boldsymbol{p}|-|\boldsymbol{p}^{\prime}|-f_{N}(\boldsymbol{p}^{\prime})\}.
\end{equation}
Here, $|\boldsymbol{p}|-|\boldsymbol{p}^{\prime}|$ means the number of stones taken in the move from $\boldsymbol{p}$ to $\boldsymbol{p}^{\prime}$, with each stone contributing one point. This recursion formula allows us to determine $f_{N}(\boldsymbol{p})$ sequentially, starting from positions with a smaller total number of stones.

\begin{figure}[h]
\centering
\hspace*{1cm}
\unitlength.15pt
\begin{picture}(1300,750)(100,100)
\put(0,200){\line(1,0){1350}}
\put(650,100){\line(0,1){700}}
\put(0,550){\line(1,-1){250}}
\put(250,300){\line(1,1){200}}
\put(450,500){\line(1,-1){100}}
\put(550,400){\line(1,1){100}}
\put(650,500){\line(1,-1){100}}
\put(750,400){\line(1,1){100}}
\put(850,500){\line(1,-1){200}}
\put(1050,300){\line(1,1){250}}
\put(600,130){\large{$0$}}
\put(600,450){\large{$3$}}
\put(180,130){\large{$-4$}}
\put(1030,130){\large{$4$}}
\put(1320,120){\large{$N$}}
\end{picture}
\caption{The plot of the function $f_{N}(5,4,2)$}
\label{fig2}
\end{figure}
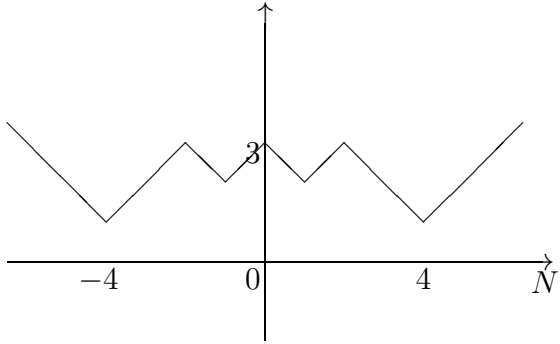

The payoff function $f_{N}(\boldsymbol{p})$ is a function of the position $\boldsymbol{p}$, but if we fix $\boldsymbol{p}$ and regard the bonus point $N$ as a variable, it can also be considered a function of $N$ with the real numbers as its domain. 
To illustrate the payoff function, Fig. \ref{fig2} shows $f_{N}(5,4,2)$ as a function of $N$ for the game that was used in the previous section. A similar position, $f_{N}(5,4,3)$, will be presented later in Remark \ref{rem1} and Fig. \ref{fig6} in Section \ref{sec6}. 
\vspace*{-1mm}

The payoff function $f_{N}(\boldsymbol{p})$ satisfies the following properties. 

\medskip

\begin{proposition}\label{prop1} 
(i) If $N$ is an integer, then $f_{N}(\boldsymbol{p})$ is also an integer, and its parity coincides with that of $|\boldsymbol{p}|+N$. \\
(ii) When $f_{N}(\boldsymbol{p})$ is regarded as a function of $N$, 
it is continuous in $N$, and for non-integer values of $N$, its slope is either $1$ or $-1$. 
In the former case, the first player takes the last stone under optimal play, 
whereas in the latter case, the second player does.
\end{proposition}

\medskip

{\bf Proof.} (i) If $N$ is an integer, then the scores of the first and second players are integers, and their sum is $|\boldsymbol{p}|+N$. 
Observe that
\[
f_{N}(\boldsymbol{p})=|\boldsymbol{p}|+N-2(\mbox{Second player's score}),
\]
and hence the parity coincides with that of $|\boldsymbol{p}|+N$. \\
(ii) We prove the claims by induction on the total number of stones  $|\boldsymbol{p}|$. 
First consider  $|\boldsymbol{p}|=0$, that is, the terminal position. In this case $f_N(\boldsymbol{0})=-N$, which is a continuous function of $N$ with slope $-1$.
Since this position can be regarded as the one obtained when the previous player takes the last stone and the game ends, assertion (ii) holds.

Assume that the statements hold for all positions with $|\boldsymbol{p}|<m$.
For $|\boldsymbol{p}|=m$, the recursion formula (\ref{eq2}) shows that $f_N(\boldsymbol{p})$ is defined as the maximum of finitely many continuous functions, each having slope $\pm 1$. Hence it is continuous and its slope is either $1$ or $-1$.
Moreover, for a non-integer $N$, whenever there exists a position $\boldsymbol{p'}$ such that $f_N(\boldsymbol{p})=|\boldsymbol{p}|-|\boldsymbol{p'}|-f_{N}(\boldsymbol{p'})$, the slope of $f_N(\boldsymbol{p})$ has the opposite sign to that of $f_N(\boldsymbol{p'})$. This corresponds to the alternation of turns, and therefore the player who takes the last stone is determined according to this sign. \quad $\Box$

\medskip

By this proposition, when considering $f_{N}(\boldsymbol{p})$ as a function of $N$, if the values for integer $N$ are determined, then the values for all real $N$ can also be obtained by linear interpolation. 

The following fact follows from Proposition \ref{prop1} (ii), since it implies that the absolute value of the slope of  $f_N(\boldsymbol{p})$ is at most 1. As this fact is useful, we present it here as a corollary.

\medskip

\begin{corollary}\label{cor2}
For real numbers $N$ and $N'$, $|f_{N}(\boldsymbol{p})-f_{N'}(\boldsymbol{p})|\leq |N-N'|$ holds.
\end{corollary}

\section{Scoring Nim with two piles}\label{sec4}
In this section, we derive an explicit formula for $f_{N}(\boldsymbol{p})$ when the initial position has two piles. Using equations (\ref{eq1}) and (\ref{eq2}), we first obtain the following results. 

\medskip

\begin{proposition}\label{prop3} $f_{N}(1,1)=-N$, and 
\[
\quad f_{N}(x,1)= x-1+|N|
\]
for $x \geq 2$.
\end{proposition}

\medskip

{\bf Proof.} In the position $(1,1)$, each player takes one stone; the first player earns 1 point, while the second player earns $1+N$ points. Moreover, applying (\ref{eq1}) and the recursion formula  (\ref{eq2}),
\begin{eqnarray*}
&&f_{N}(x,1)\\
&=&\max\left\{\begin{array}{l}
x-f_N(0,1),x-1-f_N(1,1),x-2-f_N(2,1),\ldots,\\1-f_N(x-1,1),
1-f_N(x,0)
\end{array}
\right\}\\
&=&\max\left\{\begin{array}{l}
x-1-N, x-1+N,x-3-|N|,\ldots,\\
-x+3+|N|,-x+2-|1+N|\end{array}
\right\}\\
&=&x-1+|N|
\end{eqnarray*}
holds for $x\geq 2$ inductively. \quad $\Box$

\medskip

In the position $(x,1)$ with $x\geq 2$, the optimal strategy depends on the value of $N$.  If $N\geq 0$, the optimal move is to take $x-1$ stones, leading to the position $(1,1)$. On the other hand, 
if $N\leq 0$, the optimal move is to take $x$ stones, leading to the position $(0,1)$. Since the number of possible moves is relatively small for the  two-pile case, $f_{N}(x,y)$ can be computed directly using the recursion formula. The following proposition completes the analysis for all cases. 

\medskip

\begin{proposition}\label{prop4} Let $x,y\geq 2$. Then $f_{N}(x,y)$ satisfies the following.\\
(i) 
\begin{equation}\label{eq3}
f_{N}(x,x)=x-f_{N}(0,x)=1-|1+N|.
\end{equation}
(ii) If $x>y$,
\begin{equation}\label{eq4}
f_{N}(x,y)=\max\{x-f_{N}(0,y), x-y-f_{N}(y,y)\}= x-y +|1-|1+N||.
\end{equation}
\end{proposition}

\medskip

If both piles contain the same number of stones, the optimal move is to take all stones from one pile.
When the piles have different sizes and $-2\leq N\leq 0$, the optimal move follows the greedy strategy of taking all stones from the larger pile.
Otherwise, the optimal move follows the Nim-winning strategy of equalizing the numbers of stones in the two piles.
The proof follows the same approach as in Proposition \ref{prop3} and can be established by induction on the total number of stones, using the recursion formula. 

\medskip

{\bf Proof. } We prove (i) and (ii) for all $(x,y)$ satisfying $x\geq y\geq 2$ by induction. 
Assume that (\ref{eq3}) and (\ref{eq4}) hold for every $(x',y')$ with $2\leq x'\leq x$, $2\leq y'\leq y$ and $x'+y'<x+y$.
Then, together with (\ref{eq1}) and Proposition \ref{prop3}, it follows that the following equalities hold:
\begin{eqnarray*}
&&f_N(x,x)\\
&=&\max\{x-f_N(x,0),x-1-f_N(x,1),x-2-f_N(x,2),\ldots,1-f_N(x,x-1)\}\\
&=&\max\{1-|1+N|,-|N|,-|1-|1+N||,\ldots,-|1-|1+N||\}\\
&=&1-|1+N|,\\
&&f_N(x,y)\\
&=&\max\left\{\begin{array}{l}
x-f_N(0,y),x-1-f_N(1,y),\ldots,x-y-f_N(y,y),\ldots,\\
1-f_N(x-1,y),\\
y-f_N(x,0),y-1-f_N(x,1),y-2-f_N(x,2),\ldots,\\
1-f_N(x,y-1)
\end{array}
\right\}\\
&=&\max\left\{\begin{array}{l}
x-y+1-|1+N|,x-y-|N|,\ldots,x-y-1+|1+N|,\ldots, \\
-x+y+2-|1-|1+N||,\\
-x+y+1-|1+N|,-x+y-|N|,-x+y +|1-|1+N||,\ldots,\\
-x+y +|1-|1+N||
\end{array}
\right\}\\
&=&x-y +|1-|1+N||. 
\end{eqnarray*}
In the last maximization step, we use the assumption that $x>y$. \quad $\Box$

\medskip

These results indicate that Scoring Nim with two piles has well-defined strategic properties.

\section{General properties of the payoff function}\label{sec5}
When there are three or more piles, explicitly determining $f_{N}(\boldsymbol{p})$ becomes challenging. In this section, we establish several propositions regarding $f_{N}(\boldsymbol{p})$ in Scoring Nim. The following two propositions are obtained by dividing the piles into two groups and analyzing the game as a disjunctive sum. 

\medskip

\begin{proposition}\label{prop5} 
\[f_{N}(p_{1},p_{2},\ldots,p_{n},1,1)=f_{N}(p_{1},p_{2},\ldots,p_{n}).\]
\end{proposition}

This proposition simplifies the analysis by removing pairs of piles containing only one stone. 

\medskip

{\bf Proof. } Let $\boldsymbol{p}=(p_{1},p_{2},\ldots,p_{n},1,1)$. Divide the piles into $\boldsymbol{\overline{p}}=(p_{1},p_{2},\ldots,p_{n})$ and $\boldsymbol{\underline{p}}=(1,1)$, and let $\alpha=f_{N}(\boldsymbol{\overline{p}})$.

We first show that the second player can ensure $f_{N}(\boldsymbol{p})\leq \alpha$. If the first player makes a move in $\boldsymbol{\overline{p}}$, the second player responds optimally within $\boldsymbol{\overline{p}}$. If the first player makes a move in $\boldsymbol{\underline{p}}$,  the second player removes the remaining stone from $\boldsymbol{\underline{p}}$. Furthermore, if the first player takes the last stone from $\boldsymbol{\overline{p}}$, the second player removes a remaining stone from $\boldsymbol{\underline{p}}$, if possible. In this case, the first player will then remove the remaining stone from $\boldsymbol{\underline{p}}$, ending the game. As a result, the stones in $\boldsymbol{\underline{p}}$ are removed one by one by each player, meaning that the identity of the player who takes the last stone in $\boldsymbol{\overline{p}}$ remains the same in $\boldsymbol{p}$. Therefore, the second player can ensure that $f_{N}(\boldsymbol{p})\leq \alpha$.

Similarly, we show that the first player can ensure $f_{N}(\boldsymbol{p})\geq \alpha$. The first player begins by making the optimal move in $\boldsymbol{\overline{p}}$. If the second player makes a move in $\boldsymbol{\overline{p}}$, the first player responds optimally within $\boldsymbol{\overline{p}}$ if there are stones remaining in $\boldsymbol{\overline{p}}$. If the second player makes a move in $\boldsymbol{\underline{p}}$, the first player removes the remaining stone from $\boldsymbol{\underline{p}}$. By following this strategy, the first player ensures $f_{N}(\boldsymbol{p})\geq \alpha$. \quad $\Box$

\medskip

The following proposition states that if a position is symmetric, meaning it consists of two identical groups of piles, the same payoff function as (\ref{eq3}) is obtained.

\medskip

\begin{proposition}\label{prop6} Let the number of piles $n$ be even (i.e., $n=2m$), and suppose that in $\boldsymbol{p}=(p_{1},p_{2},\ldots,p_{n})$, the condition $p_{2j-1}=p_{2j}$ holds for $1\leq j\leq m$. If there exists $i$ such that $p_{i}\geq 2$,
\[f_{N}(\boldsymbol{p})=1-|1+N|\]
holds. 
\end{proposition}

\medskip

{\bf Proof.} By Proposition \ref{prop5}, we may assume that all $p_{i}$ are at least 2. For $m=1$, the result follows directly from  (\ref{eq3}).   We prove the claim for $m\geq 2$ by induction on $|\boldsymbol{p}|$. Let $\alpha=1-|1+N|$. 

First, we establish a copy-cat strategy for the second player to ensure $f_{N}(\boldsymbol{p})\leq \alpha$. 
When the first player takes stones from a pile containing $p_{i}$ stones, the second player responds by taking the same number of stones from another pile that also originally contained $p_{i}$ stones. This restores the symmetry of the position, and by the induction hypothesis, we conclude that $f_{N}(\boldsymbol{p})\leq \alpha$. 

Next, we construct a strategy for the first player to ensure $f_{N}(\boldsymbol{p})\geq \alpha$. The first player begins by taking all stones from the pile with the largest number of stones. If the second player then takes the same number of stones from another pile, symmetry is restored. If instead, the second player takes fewer stones than the first player on the initial move, the first player can continue to always take all stones from the pile with the largest remaining number of stones. By the end of the game, the total number of stones taken by the first player will be at least two more than that taken by the second player, ensuring that the final score of the first player is at least $2-|N|$, which is at least $\alpha$. \quad $\Box$

\medskip

When $N$ is 0, the greedy strategy is optimal. By slightly extending the range of $N$, we obtain the following proposition. 

\medskip

\begin{proposition}\label{prop7} 
Let $\boldsymbol{p}=(p_{1},p_{2},\ldots,p_{n})$ with $p_{1}\geq p_{2}\geq\cdots\geq p_{n}>0$. Then, for $-1\leq N\leq 0$, the following  hold: \\
(i) When $n$ is odd, i.e., $n=2k-1$, we have
\[
f_{N}(\boldsymbol{p})=\sum_{j=1}^{k-1}\left(p_{2j-1}-p_{2j}\right)+p_{2k-1}+N.
\]
(ii) When $n$ is even, i.e., $n=2k$, we have
\[
f_{N}(\boldsymbol{p})=\sum_{j=1}^{k}\left(p_{2j-1}-p_{2j}\right)-N.
\]
\end{proposition}

\medskip

{\bf Proof.} Scoring Nim for $-1\leq N\leq 0$ can be viewed as a variant of the stone-taking game in which only the last stone is worth $1+N$ points  (with $0\leq 1+N\leq 1$) instead of 1.  Since each player aims to take as many stones as possible, this proposition follows almost trivially. \quad $\Box$

\medskip

Next, we consider the case when $|N|$ is very large. If $N$ is positive, both players will attempt to follow the normal play Nim strategies to take the last stone. Conversely, if $N$ is negative, they will attempt to follow the mis\`ere play Nim strategies to force the opponent to take the last stone. To elaborate on this, we first revisit the sets of P-positions in both the normal and mis\`ere play versions of Nim. For non-negative integers $a$ and $b$, we express them in binary notation as 
$a=\sum_{i} 2^{i}a_{i}, b=\sum_{i}2^{i}b_{i}$. The Nim-sum (exclusive OR) of $a$ and $b$, denoted as $a\oplus b$, is defined as $\displaystyle{a\oplus b=\sum_{i}2^{i}((a_{i}+b_{i})\! \! \mod 2)}$, which is the sum without carry in binary representation. Since the identity $(a\oplus b)\oplus c=a\oplus (b\oplus c)$ holds, the Nim-sum can be naturally extended to three or more non-negative integers. 

\medskip

\begin{definition}[P-position set in Nim]\label{def2} Define the sets of positions $\boldsymbol{p}=$\\
$(p_{1},p_{2},\ldots,p_{n})$ as follows:\begin{eqnarray*}
P^{+}&=& \{ \ \boldsymbol{p} \ | \ p_{1}\oplus p_{2}\oplus\cdots\oplus p_{n}=0 \ \},\\
P^{-}&=& \{ \ \boldsymbol{p} \ | \ \mbox{there exists } i \mbox{ such that } p_{i}\geq {2}, \mbox{ and } p_{1}\oplus p_{2}\oplus\cdots\oplus p_{n}=0 \ \}\\
&&\cup \ \{ \ \boldsymbol{p} \ | \ p_{i}\leq 1 \mbox{ for all } i, \mbox{ and } |\boldsymbol{p}|\mbox{ is odd} \ \}.
\end{eqnarray*}
Here, $P^{+}$ represents the set of P-positions in the normal play Nim, while $P^{-}$ represents the set of P-positions in the mis\`ere play Nim. 
\end{definition}

\medskip

Using these P-position sets, we describe the behavior of the payoff function when $|N|$ is sufficiently large.

\medskip

\begin{proposition}\label{prop8} 
(i) There exists an integer $c^{+}(\boldsymbol{p})$ such that
\[
f_{N}(\boldsymbol{p})=
\left\{
\begin{array}{ll}
c^{+}(\boldsymbol{p})-N & \mbox{ if } \boldsymbol{p}\in P^{+}, \\
c^{+}(\boldsymbol{p})+N & \mbox{ if } \boldsymbol{p}\not\in P^{+}
\end{array}
\right.
\]
holds for $N\geq |\boldsymbol{p}|-2$.\vspace{2mm}\\
(ii) There exists an integer $c^{-}(\boldsymbol{p})$ such that 
\[
f_{N}(\boldsymbol{p})=
\left\{
\begin{array}{ll}
c^{-}(\boldsymbol{p})+N & \mbox{ if } \boldsymbol{p}\in P^{-}, \\
c^{-}(\boldsymbol{p})-N & \mbox{ if } \boldsymbol{p}\not\in P^{-}
\end{array}
\right.
\]
holds for $N\leq -|\boldsymbol{p}|$.
\end{proposition}

\medskip

This proposition states that, for example, when $N$ is positive and sufficiently large, the player who takes the last stone is determined by whether the position is in $P^{+}$. Under the constraint that the player taking the last stone must always move to a P-position, the difference in the number of stones taken by both players assuming they take as many stones as possible is given by $c^{+}(\boldsymbol{p})$. For instance, in the position $(5,4,2)$, if $N$  is positive and sufficiently large, the first player will take the last stone. However, to achieve this, the first player must allow the second player to take 7 stones, which is unavoidable (while first player takes 4 stones). Thus, the payoff function is given by $f_{N}(5,4,2)=(4-7)+N=-3+N$,  implying that $c^{+}(5,4,2)=-3$.

\medskip

{\bf Proof of Proposition \ref{prop8}.} (i) When $\boldsymbol{p}\not\in P^{+}$, the first player has a winning strategy in the normal play Nim (i.e., taking the last stone), ensuring they obtain at least $1+N$ points. On the other hand, if the first player does not take the last stone, the total points obtained by taking all remaining stones is at most $|\boldsymbol{p}|-1$. Thus, if $1+N\geq |\boldsymbol{p}|-1$, which simplifies to $N\geq |\boldsymbol{p}|-2$, the first player will always select to take the last stone. Under the constraint that the first player must always select the next game state from within $P^{+}$, an optimal strategy that maximizes their payoff is determined. The value of $c^{+}(\boldsymbol{p})$ is given by the difference between the number of stones the first player can take and the number the second player can take. If $\boldsymbol{p}\in P^{+}$, then $c^{+}(\boldsymbol{p})$ is determined under the constraint that the second player must always select the next game state from within $P^{+}$. \\
(ii) In this case, similar to (i), we compare the maximum possible points when a player takes the last stone, $|\boldsymbol{p}|+N$, with the minimum possible points when a player avoids to take it, which is 0. If  $|\boldsymbol{p}|+N\leq 0$, the player who wins in the  mis\`ere play Nim will always select the next game state from within $P^{-}$. Under this constraint,  $c^{-}(\boldsymbol{p})$ is determined. \quad $\Box$\vspace{2mm}

As an example of $c^{+}(\boldsymbol{p})$ and $c^{-}(\boldsymbol{p})$, we compute $c^{+}(3,2,1)$ and $c^{-}(3,2,1)$ in the three-pile Scoring Nim. Starting from $(3,2,1)$, the opponent responds with a move that  leads to the next game state in either $P^{+}$ or $P^{-}$ position when $|N|$ is sufficiently large. Specifically, possible moves include $(3,2,1)\to (2,2,1)\to(2,2,0)\in P^{+}\cap P^{-}$, $(3,2,1)\to(1,2,1)\to(1,0,1)\in P^{+}$ or $(1,1,1)\in P^{-}$, among others.  By evaluating all possible moves, we find that:
\[
c^{+}(3,2,1)=3-1+c^{+}(0,1,1)=2, \quad c^{-}(3,2,1)=3-2+c^{-}(0,0,1)=2.
\]

Similarly, considering all possible outcomes for $(2k+1,2k,1)$ ($k\geq 2$) when $|N|$ is sufficiently large, we find that the optimal sequence of moves consists of the current player first taking 3 stones from the pile with $2k+1$ stones,  followed by the opponent taking 1 stone from the pile with $2k$ stones; that is,
\begin{eqnarray*}
c^{+}(2k+1,2k,1)&=&3-1+c^{+}(2k-2,2k-1,1), \\
c^{-}(2k+1,2k,1)&=&3-1+c^{-}(2k-2,2k-1,1).
\end{eqnarray*}
This leads to the inductive conclusion that
\begin{equation}\label{eq5}
c^{+}(2k+1,2k,1)=c^{-}(2k+1,2k,1)=2k \quad (k\geq 1).
\end{equation}
Here the current player effectively implements a strategy in which, through alternating moves, they repeat as many times as possible the situation where they take two more stones than the opponent.
This result will be used in the proof of Theorem \ref{thm11} in the next section.

\section{Changes in the optimal strategy depending on the bonus value}\label{sec6}
As discussed in the previous section, when $N$ is close to 0 or when $|N|$ is sufficiently large, the problem can be effectively analyzed by focusing solely on the greedy strategy or the Nim-winning strategy, which provides a partial solution. However, for intermediate values of $N$, the game exhibits significant complexity in optimal move selection. Therefore, in this section, we  focus our analysis on three-pile positions and investigate the structure of the payoff function $f_{N}(x,y,z)$ with respect to $N$. First, we consider cases where multiple piles contain the same number of stones, which simplifies the determination of the payoff function.

\medskip

\begin{proposition}\label{prop9} The function $f_{N}(x,y,z)$ is given as follows for $x,y,z\geq 1$ when at least two of them are equal:
$f_{N}(1,1,1)=1+N$, and if at least one of $x$ or $z$ is at least 2, then 
\begin{equation}\label{eq6}
f_{N}(x,x,z)=z-1+|1+N|
\end{equation}
holds. 
\end{proposition}

\medskip

{\bf Proof.} 
For  (\ref{eq6}), we apply the recursion formula  (\ref{eq2}) and use induction on the total number of stones to establish the result. We assume here $x>z$. 

First, the first player may take all the stones from the pile of size $z$. By  (\ref{eq3}), we have
\[
f_{N}(x,x,z)\geq z-f_{N}(x,x,0)=z-(1-|1+N|).
\] 

To prove the reverse inequality, we show that for any move by the first player, the second player can respond so that $f_{N}(x,x,z)$ does not exceed $z-1+|1+N|.$
\begin{itemize}
\item If the first player moves to $(x,0,z)$, then by  (\ref{eq4}), we have
\[
x-f_{N}(x,0,z)=x-\left((x-z)+|1-|1+N||\right)\leq z-(1-|1+N|).
\]
\item If the first player moves to $(x,x',z)$ (where $0<x'<x$), the second player responds by moving to $(x',x',z)$ unless $x'=1$ and $z=1$. By the induction hypothesis, it follows that
\[x-x'-f_{N}(x,x',z) \leq x-x'-\left(x-x'-f_{N}(x',x',z)\right) =z-1+|1+N|.\]
If $x'=1$ and $z=1$, the second player responds by moving to $(0,1,1)$, and hence
\[x-1-f_{N}(x,1,1) \leq x-1-\left(x-f_{N}(0,1,1)\right) =-1-N\leq |1+N|.\] 
\item If the first player moves to $(x,x,z')$ (where $0<z'<z$), the second player responds by moving to $(x,x,0)$. Then, by  (\ref{eq3}), we have
\begin{eqnarray*}
z-z'-f_{N}(x,x,z')
&=&z-z'-\left(z'-f_{N}(x,x,0)\right)\\
&=&z-2z'+(1-|1+N|)\\
&\leq&z-1+|1+N|.
\end{eqnarray*}
\end{itemize}

Thus, the case $x>z$ has been proved.
The cases $x=z$ and $x<z$ follow since the optimal play for the first player is also to move to $(x,x,0)$.
\quad $\Box$

\medskip

This proposition narrows our focus to the payoff function $f_{N}(x,y,z)$ for $x>y>z>0$. However, fully determining this function in general cases is challenging at this stage. Therefore, we focus, in particular, on the case $z=1$. We analyze $f_{N}(x,y,1)$ and show that this game reveals intriguing strategic behavior, where the optimal strategy varies significantly depending on the value of the bonus.

We classify $f_{N}(x,y,1)$ according to cases, where the condition $(x,y,1)\in P^{+}\cap P^{-}$  holds only when $(2k+1,2k,1)$ with $k\geq 1$. We define the function $F_{k}(N)$ for $k\geq 1$ as follows.

\medskip

\begin{definition}\label{def3} 
For any integer $k\geq 1$, we define $J_{k}=\{-(2k-2),-(2k-4), \ldots, -4,-2,0,2,4,\ldots, 2k-4,2k-2\}$. Based on $J_{k}$, we define
\[
F_{k}(N)=2-\min_{j\in J_{k}} |N-j|.
\] 
\end{definition}

\medskip

\begin{figure}[h]
\centering
\hspace*{1cm}
\unitlength.15pt
\begin{picture}(1300,500)(100,200)
\put(0,300){\line(1,0){1350}}
\put(650,200){\line(0,1){500}}
\put(0,250){\line(1,1){250}}
\put(250,500){\line(1,-1){100}}
\put(350,400){\line(1,1){100}}
\put(450,500){\line(1,-1){100}}
\put(550,400){\line(1,1){100}}
\put(650,500){\line(1,-1){100}}
\put(750,400){\line(1,1){100}}
\put(850,500){\line(1,-1){100}}
\put(950,400){\line(1,1){100}}
\put(1050,500){\line(1,-1){250}}
\put(600,230){\large{$0$}}
\put(0,230){\large{$-6$}}
\put(1230,230){\large{$6$}}
\put(600,450){\large{$2$}}
\put(1320,220){\large{$N$}}
\end{picture}
\caption{The plot of the function $F_{3}(N)$}
\label{fig3}
\end{figure}
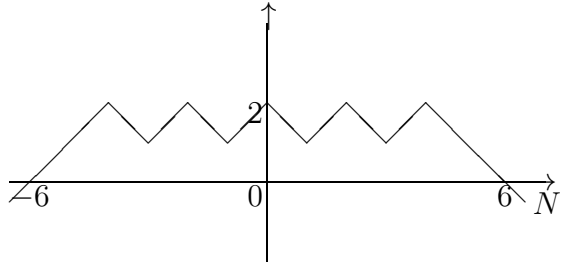

For example, when $k=3$, $F_{3}(N)$  is depicted in Fig. \ref{fig3}.  We record the following properties of $F_{k}(N)$ as a lemma.

\medskip

\begin{lemma}\label{lem10}  For $k\geq 1$, the function $F_{k}(N)$ satisfies the following properties:\\
(i) $F_{k}(N)\leq 2$ for all $N$, and $F_{k}(N)=2$ if and only if $N\in J_{k}$.\\
(ii) When $N$ is an odd integer with $|N|\leq 2k-1$, then $F_{k}(N)=1$.  \\ 
(iii) $F_{k}(2k)=F_{k}(-2k)=0$.\\
(iv) $F_{k}(N)=2k-|N|$ for $|N|\geq 2k-2$. \\
(v) For all $N$, $\max\{F_{k}(N),-2k+|N| \}\geq 0$ and $\max\{F_{k}(N),-(2k-2)+|N|\}\geq 1$.
\end{lemma}

\medskip

Based on $F_{k}(N)$, $f_{N}(x,y,1)$ can be formulated in terms of $N$ as follows.

\medskip

\begin{theorem}\label{thm11} (i) For any integer $k\geq 1$,
\[
f_{N}(2k+1,2k,1)=F_{k}(N).
\]
(ii) For any integers $k\geq 1$ and $x\geq 2k+2$,
\[
f_{N}(x,2k,1)=\max\{F_{k}(N),-2k+|N| \}+x-2k-1.
\]
(iii) For any integers $k\geq 1$ and $x\geq 2k+2$, 
\[
f_{N}(x,2k+1,1)=\max\{F_{k}(N),-(2k-2)+|N|\}+x-2k-2.
\]
\end{theorem}

\medskip

\begin{figure}[h]
\centering
\hspace*{1cm}
\unitlength.15pt
\begin{picture}(1300,600)(100,150)
\put(0,300){\line(1,0){1350}}
\put(650,200){\line(0,1){500}}

\put(-5,235){$\cdot$}
\put(5,245){$\cdot$}
\put(15,255){$\cdot$}
\put(25,265){$\cdot$}
\put(35,275){$\cdot$}
\put(45,285){$\cdot$}
\put(55,295){$\cdot$}
\put(65,305){$\cdot$}
\put(75,315){$\cdot$}
\put(85,325){$\cdot$}
\put(95,335){$\cdot$}
\put(105,345){$\cdot$}
\put(115,355){$\cdot$}
\put(125,365){$\cdot$}
\put(135,375){$\cdot$}

\put(145,375){$\cdot$}
\put(155,365){$\cdot$}
\put(165,355){$\cdot$}
\put(175,345){$\cdot$}
\put(185,335){$\cdot$}
\put(195,325){$\cdot$}
\put(205,315){$\cdot$}
\put(215,305){$\cdot$}
\put(225,295){$\cdot$}
\put(235,285){$\cdot$}
\put(245,275){$\cdot$}
\put(255,265){$\cdot$}
\put(265,255){$\cdot$}
\put(275,245){$\cdot$}
\put(285,235){$\cdot$}
\put(295,225){$\cdot$}
\put(305,215){$\cdot$}
\put(315,205){$\cdot$}
\put(325,195){$\cdot$}
\put(335,185){$\cdot$}

\put(945,185){$\cdot$}
\put(955,195){$\cdot$}
\put(965,205){$\cdot$}
\put(975,215){$\cdot$}
\put(985,225){$\cdot$}
\put(995,235){$\cdot$}
\put(1005,245){$\cdot$}
\put(1015,255){$\cdot$}
\put(1025,265){$\cdot$}
\put(1035,275){$\cdot$}
\put(1045,285){$\cdot$}
\put(1055,295){$\cdot$}
\put(1065,305){$\cdot$}
\put(1075,315){$\cdot$}
\put(1085,325){$\cdot$}
\put(1095,335){$\cdot$}
\put(1105,345){$\cdot$}
\put(1115,355){$\cdot$}
\put(1125,365){$\cdot$}
\put(1135,375){$\cdot$}

\put(1145,375){$\cdot$}
\put(1155,365){$\cdot$}
\put(1165,355){$\cdot$}
\put(1175,345){$\cdot$}
\put(1185,335){$\cdot$}
\put(1195,325){$\cdot$}
\put(1205,315){$\cdot$}
\put(1215,305){$\cdot$}
\put(1225,295){$\cdot$}
\put(1235,285){$\cdot$}
\put(1245,275){$\cdot$}
\put(1255,265){$\cdot$}
\put(1265,255){$\cdot$}
\put(1275,245){$\cdot$}
\put(1285,235){$\cdot$}

\put(600,230){\large{$0$}}
\put(0,230){\large{$-6$}}
\put(200,230){\large{$-4$}
\put(730,0){\large{$4$}}
\put(930,0){\large{$6$}}}
\put(600,450){\large{$2$}}
\put(1320,220){\large{$N$}}
\put(0,550){\line(1,-1){150}}
\put(150,400){\line(1,1){100}}
\put(250,500){\line(1,-1){100}}
\put(350,400){\line(1,1){100}}
\put(450,500){\line(1,-1){100}}
\put(550,400){\line(1,1){100}}
\put(650,500){\line(1,-1){100}}
\put(750,400){\line(1,1){100}}
\put(850,500){\line(1,-1){100}}
\put(950,400){\line(1,1){100}}
\put(1050,500){\line(1,-1){100}}
\put(1150,400){\line(1,1){150}}
\end{picture}
\caption{$f_{N}(8,7,1)=\max\{F_{3}(N),-4+|N|\}$ (solid line)}
\label{fig4}
\end{figure}
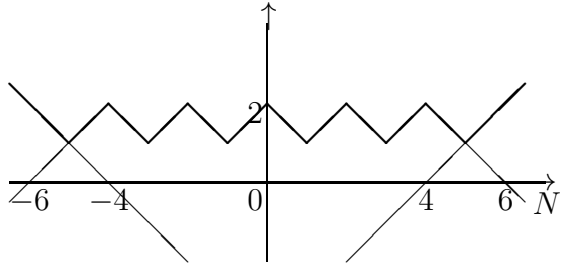

Thus, we have fully characterized $f_{N}(x,y,1)$. For example, $f_{N}(8,7,1)$ is an instance of the case where $k=3$ and $x=8$ in part (iii) of the theorem above. It is given by $\max\{F_{3}(N),-4+|N|\}$, which is represented by the solid line in Fig. \ref{fig4}. 

One key observation concerns the number of breakpoints in $f_{N}(x,y,1)$, namely, points where the slope changes sign. We find that $f_{N}(2k+1,2k,1)$ has $4k-3$ breakpoints, indicating that the number of breakpoints can grow arbitrarily large as the number of stones increases.

\medskip

{\bf Proof of Theorem \ref{thm11}.} For $k=1$, (i), (ii), and (iii) can be verified directly. For $k\geq 2$, we proceed by induction, assuming that (i), (ii), and (iii) hold for all $k<m$, and proving that (i) holds for $k=m$. Under this assumption, the following facts (a), (b), and (c) hold for $1\leq k <m$.\vspace{6pt}\\
(a) For all $N$, we have $f_{N}(2k+1,2k,1)\leq 2$, with equality holding for all $N$ in $J_{k}$. In particular, if $N$ is an odd integer, then $f_{N}(2k+1,2k,1)\leq 1$. Moreover, if $|N|\geq 2k-2$, then $f_{N}(2k+1,2k,1)=2k-|N|$. Consequently, $f_{N}(2k+1,2k,1)\leq -1$ if $|N|\geq 2k+1$ and $f_{N}(2k+1,2k,1)\leq -2$ if $|N|\geq 2k+2$.\\
(b) For $x\geq 2k+2$, it follows that $f_{N}(x,2k,1)\geq x-2k-1$ for all $N$, with equality holding only when $N=\pm 2k$. Moreover, when $N=\pm(2k+1)$, we have $f_{N}(x,2k,1)=x-2k$.\\
(c) For $x\geq 2k+2$, it follows that $f_{N}(x,2k+1,1)\geq x-2k-1$ for all $N$.

\medskip

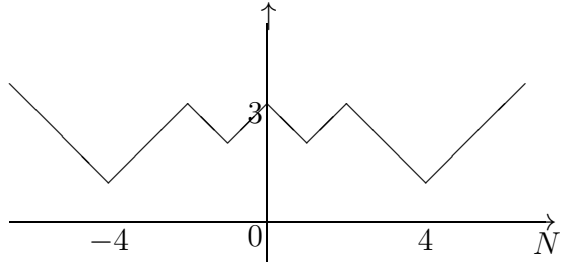
\begin{figure}[h]
\centering
\hspace*{1cm}
\unitlength.15pt
\begin{picture}(1300,650)(100,100)
\put(0,200){\line(1,0){1350}}
\put(650,100){\line(0,1){600}}
\put(0,550){\line(1,-1){250}}
\put(250,300){\line(1,1){200}}
\put(450,500){\line(1,-1){100}}
\put(550,400){\line(1,1){100}}
\put(650,500){\line(1,-1){100}}
\put(750,400){\line(1,1){100}}
\put(850,500){\line(1,-1){200}}
\put(1050,300){\line(1,1){250}}
\put(600,130){\large{$0$}}
\put(200,130){\large{$-4$}}
\put(1030,130){\large{$4$}}
\put(600,450){\large{$3$}}
\put(1320,120){\large{$N$}}
\end{picture}
\caption{The plot of the function $f_{N}(6,4,1)$}
\label{fig5}
\end{figure}

These facts follow easily from Lemma \ref{lem10} and will be used repeatedly in this proof. Notably, in (b), equality holds for only two values of $N$. Specifically, $f_{N}(6,4,1)$ is an instance of the case where $k=2$ and $x=6$ in Theorem \ref{thm11} (ii), and as illustrated in Fig. \ref{fig5}, it attains its minimum value at $N=\pm 4$.

As stated in Section 3, Proposition 1(ii) implies that $f_{N}(p)$ is a continuous piecewise linear function of $N$ with slope $\pm 1$.
Hence, to establish (i), it suffices to consider integer values of $N$. We outline the proof strategy as follows:
\begin{itemize}
\item By showing that $f_{N}(2m+1,2m,1)\geq 2$ when $N$ is an even integer in the range $|N|\leq 2m-2$, it follows that $f_{N}(2m+1,2m,1)\geq F_{m}(N)$ for all $N$.
\item By showing that $f_{N}(2m+1,2m,1)\leq 1$ when $N$ is an odd integer in the range $|N|\leq 2m-3$, and using (\ref{eq5}), it follows that $f_{N}(2m+1,2m,1)\leq F_{m}(N)$ for all $N$.
\end{itemize}
We will explain in detail in the proof below that it is sufficient to establish these facts.

\smallskip

Case 1: $N=\pm 2l$ with $0\leq l<m$. The first player has the option to move from the position $(2m+1,2m,1)$ to $(2l,2m,1)$. By the induction hypothesis---namely, by applying Proposition \ref{prop3} when $l=0$, and using (b) when $1\leq l<m$---we have $f_{\pm 2l}(2l,2m,1)=2m-2l-1$. Hence we obtain
\[f_{\pm 2l}(2m+1,2m,1)\geq 2m+1-2l-f_{\pm 2l}(2l,2m,1)=2.
\]
Thus, it follows that $f_{N}(2m+1,2m,1)\geq 2$ for $N=-(2m-2),-(2m-4),\ldots, 2m-4,2m-2$. Furthermore, since the slope of this function with respect to $N$ is at least $-1$  (as stated in Proposition 1(ii)), it follows that $f_{N}(2m+1,2m,1)\geq F_{m}(N)$ holds for all $N$.

\smallskip

Case 2: $N=\pm(2l+1)$ with $0\leq l<m-1$. We analyze the possible moves available to the first player at the position $(2m+1,2m,1)$ and explicitly describe the opponent's responses. Through this analysis, we show that $f_{\pm(2l+1)}(2m+1,2m,1)$ does not exceed 1.  
\begin{itemize}
\item If the first player moves to $(2m+1,2m,0)$, then by (\ref{eq4})---with the second player responding by moving to $(2m,2m,0)$---we obtain
\[1-f_{\pm(2l+1)}(2m+1,2m,0)=f_{\pm(2l+1)}(2m,2m,0)=1-|1+N|\leq 1.\]
Similarly, if the first player moves to $(2m,2m,1)$, the second player responds by moving to $(2m,2m,0)$.
\item If the first player moves to $(0,2m,1)$, the second player responds by moving to $(0,1,1)$, which is optimal, when $N=2l+1$. Then
\[
2m+1-f_{2l+1}(0,2m,1)=2+f_{2l+1}(0,1,1)=2-N.
\]
If $N=-(2l+1)$, the second player responds by moving to $(0,0,1)$, which is optimal. Then
\[
2m+1-f_{-(2l+1)}(0,2m,1)=1+f_{-(2l+1)}(0,0,1)=2+N.
\]
Thus, in both cases, the function value does not exceed 1.
Similarly, if the first player moves to $(1,2m,1)$, $(2m+1,0,1)$ or $(2m+1,1,1)$, the second player responds by moving to $(1,0,1)$ if $N=2l+1$, and to either $(1,1,1)$ or $(0,0,1)$ if $N=-(2l+1)$, so that the function value does not exceed 1.
\item If  the first player moves to $(x,2m,1)$ or $(2m+1,x,1)$ (where $x$ is an even number with $2\leq x \leq 2l$), 
the second player responds by moving to $(x,x+1,1)$ (in this sequence of moves, it is more advantageous to have initially moved to $(x,2m,1)$ for the first player). Then 
\[
2m+1-x-f_{\pm(2l+1)}(x,2m,1)\leq 2+f_{\pm(2l+1)}(x,x+1,1).
\]
By the induction hypothesis (a), we have $f_{\pm(2l+1)}(x,x+1,1)\leq -1$, implying that the right-hand side of the inequality does not exceed 1.
\item If the first player moves to $(x,2m,1)$ or $(2m+1,x,1)$ (where $x$ is an odd number with $3\leq x\leq 2m-1$),  the second player responds by moving to $(x,x-1,1)$  (it is more advantageous to have initially moved to $(x,2m,1)$ for the first player). Then 
\[
2m+1-x-f_{\pm(2l+1)}(x,2m,1)\leq f_{\pm(2l+1)}(x,x-1,1).
\]
By the induction hypothesis (a), 
it follows that the right-hand side of the above inequality  does not exceed 1.
\item If the first player moves to $(x,2m,1)$ or $(2m+1,x,1)$ (where $x$ is an even number with $2l+2\leq x \leq 2m-2$), 
the second player  responds by moving to $(x,2l,1)$  (it is more advantageous to have initially moved to $(x,2m,1)$ for the first player). Then
\[2m+1-x-f_{\pm(2l+1)}(x,2m,1)\leq 1-x+2l+f_{\pm(2l+1)}(x,2l,1).\]
By the induction hypothesis (b), it follows that the right-hand side of the above inequality  does not exceed 1.
\end{itemize}

\medskip

From the above moves and responses, we have established that $f_{N}(2m+1,2m,1)\leq 1$ for $N=-(2m-3),-(2m-5),\ldots,2m-5,2m-3$. Consequently, it follows that $f_{N}(2m+1,2m,1)\leq F_{m}(N)$ for $|N|\leq 2m-2$. On the other hand, by (\ref{eq5}), it is known that for sufficiently large $|N|$, we have $f_{N}(2m+1,2m,1)=2m-|N|$. Thus, by Corollary \ref{cor2}, we conclude that for all $N$, $f_{N}(2m+1,2m,1)\leq 2m-|N|$. Hence, for all $N$, it follows that $f_{N}(2m+1,2m,1)\leq F_{m}(N)$.

Similarly, statements (ii) and (iii) for $k=m$ follow under the assumption that (i) holds for all $k\leq m$, and that (ii) and (iii) hold for all $k<m$. The details of the proof are given in the Appendix. \quad $\Box$

\medskip

We provide additional insights into Theorem \ref{thm11} and examine the factors contributing to the complexity of move selection and the increase in the number of breakpoints. First, as observed in the proof of (i), in the position $(2k+1,2k,1)$ ($k\geq 1$), the best move (or one of the best moves) for the player is to choose the even number $j$ which is closest to $N$ within the range $-(2k-2)\leq j \leq 2k-2$ and move to $(|j|,2k,1)$. For example, in the position $(9,8,1)$, the player should make the following moves based on the value of $|N|$:
\begin{itemize}
\item Moves to $(0,8,1)$ when  $0\leq |N| \leq 1$, 
\item Moves to $(2,8,1)$ when  $1\leq |N| \leq 3$,
\item Moves to $(4,8,1)$ when  $3\leq |N| \leq 5$,
\item Moves to $(6,8,1)$ when  $5\leq |N|$.   
\end{itemize}
This strategy is derived from (b), which states that $f_{N}(x,|j|,1)$ attains its minimum value $x-|j|-1$ at $N=\pm |j|$. This represents the weakest positions for the current player. 
Thus, it is advantageous to force the opponent into such a position.

Additionally, in the positions $(2k+1,2k,1)$ and $(x,2k,1)$ ($x\geq 2k+2$), considered in (i) and (ii), a common strategy is observed when $|N|$ is small: the player takes stones from the largest pile so that the remaining number of stones is close to $|N|$. As a result, the function $F_{k}(N)$, with an upward or downward shift, appears in both payoff functions. The distinction between these two payoff functions lies in whether the player can force the opponent into a P-position when $|N|$ is large. Specifically, when $|N|\geq 2k$, in the position $(2k+1,2k,1)$, the optimal move is to $(2k-2,2k,1)$, which is the same as in the case where $2k-3\leq |N|\leq 2k$. However, in the position $(x,2k,1)$ ($x\geq 2k+2$), the best move is instead to $(2k+1,2k,1)$, which is an element of $P^{+}\cap P^{-}$. For example, if the initial position is $(10,8,1)$ and $6\leq|N|\leq 8$, the game proceeds\footnote{For $5\leq |N|\leq 6$, the game  proceeds the same way up to $(10,8,1)\to(6,8,1)$, but then the second player makes a move to $(6,4,1)$.} as follows:
\[(10,8,1)\to(6,8,1)\to(6,7,1)\to(6,4,1)\to\cdots.\]
Since the first player gains 14 points and the second player gains $5+N$ points, it follows that $f_{N}(10,8,1)=9-N$. On the other hand, when 
$|N|\geq 8$, the game proceeds as follows:
\[(10,8,1)\to(9,8,1)\to(6,8,1)\to(6,7,1)\to(6,4,1)\to\cdots.\]
In this case, the first player gains $6+N$ points and the second player gains 13 points, which results in $f_{N}(10,8,1)=N-7$.

\medskip

\begin{figure}[h]
\centering
\hspace*{1cm}
\unitlength.15pt
\begin{picture}(1300,800)(100,0)
\put(0,200){\line(1,0){1350}}
\put(650,100){\line(0,1){700}}
\put(0,650){\line(1,-1){250}}
\put(250,400){\line(1,1){200}}
\put(450,600){\line(1,-1){100}}
\put(550,500){\line(1,1){100}}
\put(650,600){\line(1,-1){300}}
\put(950,300){\line(1,1){350}}
\put(600,130){\large{$0$}}
\put(600,550){\large{$4$}}
\put(180,130){\large{$-4$}}
\put(930,130){\large{$3$}}
\put(1320,120){\large{$N$}}
\end{picture}
\caption{The plot of the function $f_{N}(5,4,3)$}
\label{fig6}
\end{figure}
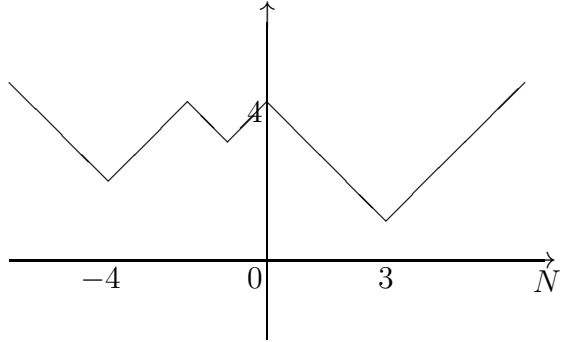

\begin{remark} \label{rem1} 
While the payoff functions presented in Theorem 11 are symmetric with respect to $N=0$, payoff functions, in general, are not symmetric. In fact, as shown in Proposition 9, $f_{N}(x,x,z)$ is symmetric with respect to $N=-1$ rather than $N=0$, and $f_{N}(5,4,3)$ (as illustrated in Fig. \ref{fig6}) lacks symmetry with respect to either $N=0$ or $N=-1$.
\end{remark}

\begin{appendices}\label{app}
\section{Proof of Theorem \ref{thm11}}
In this Appendix, we provide a detailed proof of statements (ii) and (iii) of Theorem \ref{thm11}, whose proofs were omitted in Section \ref{sec6}. Specifically, we show that statements (ii) and (iii) for $k=m$ follow under the assumption that (i) holds for all 
$k\leq m$, and that (ii) and (iii) hold for all $k<m$.

To establish (ii), it suffices to consider integer values of $N$.
It suffices to verify the following inequalities:
\begin{itemize}
\item $f_{N}(x,2m,1)\geq x-2m+1$ when $N$ is even and $|N|\leq 2m-2$,
\item $f_{N}(x,2m,1)\geq x-4m-1+|N|$ when $|N|\geq 2m$,
\item $f_{N}(x,2m,1)\leq x-2m$ when $N$ is odd and $|N|\leq 2m-3$,
\item $f_{N}(x,2m,1)\leq x-2m-1$ when $N=\pm 2m$.
\end{itemize}

Case 1: $N=\pm 2l$ with $0\leq l<m$. The first player has the option to move from the position $(x,2m,1)$ to $(2l,2m,1)$. By the induction hypothesis (b) stated at the beginning of the proof of Theorem \ref{thm11}, we obtain 
\[f_{\pm 2l}(x,2m,1)\geq x-2l-f_{\pm 2l}(2l,2m,1)=x-2m+1.\]

Case 2: $|N|\geq 2m$. The first player has the option to move from the position $(x,2m,1)$ to $(2m+1,2m,1)$. By the induction hypothesis, since (i) has already been established for $k=m$, we have $f_{N}(2m+1,2m,1)=2m-|N|$. Hence,
\[f_{N}(x,2m,1)\geq x-(2m+1)-f_{N}(2m+1,2m,1)=x-4m-1+|N|.\]

Case 3: $N=\pm(2l+1)$ with $0\leq l<m-1$. We analyze the possible moves available to the first player at the position $(x,2m,1)$ and explicitly describe the opponent's responses. Through this analysis, we show that $f_{\pm(2l+1)}(x,2m,1)$ does not exceed $x-2m$.  
\begin{itemize}
\item If the first player moves to $(2m,2m,1)$ and the second player responds by moving to $(2m,2m,0)$, then by (\ref{eq3}) we obtain
\[x-2m-f_{\pm(2l+1)}(2m,2m,1)\leq x-2m-|1+N|\leq x-2m.\]
Similarly, if the first player moves to $(x,2m,0)$, the second player responds by moving to $(2m,2m,0)$.
\item If the first player moves to $(0,2m,1)$, the second player responds by moving to $(0,1,1)$ if $N=2l+1$. Then
\[
x-f_{2l+1}(0,2m,1)\leq x-2m+1-N\leq x-2m
\]
since $N$ is a positive odd integer. 
If $N=-(2l+1)$, the second player responds by moving to $(0,0,1)$. Then
\[
x-f_{-(2l+1)}(0,2m,1)\leq x-2m+1+N\leq x-2m\]
since $N$ is a negative odd integer. 
Similarly, if the first player moves to $(1,2m,1)$, $(2m+1,0,1)$ or $(2m+1,1,1)$, the second player responds by moving to $(1,0,1)$ if $N=2l+1$, and to either $(1,1,1)$ or $(0,0,1)$ if $N=-(2l+1)$, so that the function value does not exceed $x-2m$.
\item If  the first player moves to $(x',2m,1)$ or $(x,x',1)$ (where $x'$ is an even number with $2\leq x' \leq 2l$), 
the second player responds by moving to $(x',x'+1,1)$ (in this sequence of moves, it is more advantageous to have initially moved to $(x',2m,1)$ for the first player). Then 
\[
x-x'-f_{\pm(2l+1)}(x',2m,1)\leq x-2m+1+f_{\pm(2l+1)}(x',x'+1,1).
\]
By the induction hypothesis (a), we have $f_{\pm(2l+1)}(x',x'+1,1)\leq -1$, implying that the right-hand side of the inequality does not exceed $x-2m$.
\item If the first player moves to $(x',2m,1)$ or $(x,x',1)$ (where $x'$ is an odd number with $3\leq x'\leq 2m-1$),  the second player responds by moving to $(x',x'-1,1)$  (it is more advantageous to have initially moved to $(x',2m,1)$ for the first player). Then 
\[
x-x'-f_{\pm(2l+1)}(x',2m,1)\leq x-2m-1-f_{\pm(2l+1)}(x',x'-1,1).
\]
By the induction hypothesis (a), 
it follows that the right-hand side of the above inequality  does not exceed $x-2m$.
\item If the first player moves to $(x',2m,1)$ or $(x,x',1)$ (where $2l+2\leq x'$), 
the second player  responds by moving to $(x',2l,1)$  (it is more advantageous to have initially moved to $(x',2m,1)$ for the first player). Then
\[x-x'-f_{\pm(2l+1)}(x',2m,1)\leq x-x'-2m+2l+f_{\pm(2l+1)}(x',2l,1).\]
By the induction hypothesis (b), it follows that the right-hand side of the above inequality  does not exceed $x-2m$.
\end{itemize}

Case 4: $N=\pm 2m$. A similar analysis of the available moves and responses yields $f_{\pm 2m}(x,2m,1)\leq x-2m-1$.
\begin{itemize}
\item If the first player moves to $(x,2m,0)$ or $(2m,2m,1)$, the second player responds by moving to $(2m,2m,0)$.
\item If the first player moves to $(0,2m,1)$, the second player responds by moving to $(0,1,1)$, which is optimal, when $N=2m$. 
If $N=-2m$, the second player responds by moving to $(0,0,1)$, which is optimal. 
If the first player moves to $(1,2m,1), (2m+1,0,1)$ or $(2m+1,1,1)$, the second player responds by moving to $(1,0,1)$ if $N=2m$, and to either $(1,1,1)$ or $(0,0,1)$ if $N=-2m$.  In each case, this response ensures that the function value does not exceed $x-2m-1$.
\item If  the first player moves to $(x',2m,1)$ or $(x,x',1)$ (where $x'$ is an even number with $2\leq x' \leq 2m-2$), 
the second player responds by moving to $(x',x'+1,1)$ (in this sequence of moves, it is more advantageous to have initially moved to $(x',2m,1)$ for the first player). Then 
\[
x-x'-f_{\pm 2m}(x',2m,1)\leq x-2m+1+f_{\pm 2m}(x',x'+1,1).
\]
By the induction hypothesis (a), we have $f_{\pm 2m}(x',x'+1,1)\leq -2$, implying that the right-hand side of the inequality does not exceed $x-2m-1$.
\item If the first player moves to $(x',2m,1)$ or $(x,x',1)$ (where $x'$ is an odd number with $3\leq x'\leq 2m-1$),  the second player responds by moving to $(x',x'-1,1)$  (it is more advantageous to have initially moved to $(x',2m,1)$ for the first player). Then 
\[
x-x'-f_{\pm 2m}(x',2m,1)\leq x-2m-1-f_{\pm 2m}(x',x'-1,1).
\]
By the induction hypothesis (a), 
it follows that the right-hand side of the above inequality  does not exceed $x-2m-1$.
\item If the first player moves to $(x',2m,1)$ (where $2m+2\leq x'$), 
the second player  responds by moving to $(2m+1,2m,1)$. Then
\[x-x'-f_{\pm 2m}(x',2m,1)\leq x-2x'+2m+1+f_{\pm 2m}(2m+1,2m,1).\]
By the induction hypothesis, we have $f_{\pm 2m}(2m+1,2m,0)=0$, and since $2m+2\leq x'$, the right-hand side of the above inequality  does not exceed $x-2m-1$.
\end{itemize}
\medskip

The proof of (iii) proceeds similarly.
It suffices to verify the following inequalities:

\begin{itemize}
\item $f_{N}(x,2m+1,1)\geq x-2m$ when $N$ is even and $|N|\leq 2m-2$,
\item $f_{N}(x,2m+1,1)\geq x-4m+|N|$ when $|N|\geq 2m$,
\item $f_{N}(x,2m+1,1)\leq x-2m-1$ when $N$ is odd and $|N|\leq 2m-1$.
\end{itemize}

Case 1: $N=\pm 2l$ with $0\leq l<m$. The first player has the option to move from the position $(x,2m+1,1)$ to $(2l,2m+1,1)$. By the induction hypothesis (b), we obtain 
\[f_{\pm 2l}(x,2m+1,1)\geq x-2l-f_{\pm 2l}(2l,2m+1,1)=x-2m.\]
Case 2: $|N|\geq 2m$. The first player has the option to move from the position $(x,2m+1,1)$ to $(2m,2m+1,1)$. By the induction hypothesis, since (i) has already been established for $k=m$, we obtain 
\[f_{N}(x,2m+1,1)\geq x-2m-f_{N}(2m,2m+1,1)=x-4m+|N|.\]
Case 3: $N=\pm(2l+1)$ with $0\leq l<m$. We analyze the possible moves available to the first player at the position $(x,2m+1,1)$ and explicitly describe the opponent's responses. Through this analysis, we show that $f_{\pm(2l+1)}(x,2m+1,1)$ does not exceed $x-2m-1$.  This case can be proved by the same procedure as in Case 3 of (ii).\quad$\Box$

\medskip

\end{appendices}


\medskip

\section*{Acknowledgements}
This work was supported by JSPS KAKENHI (Grant Numbers JP21K12191 and JP25K15403) and JST SPRING (Grant Number JPMJSP2115).














\end{document}